\documentclass[openany, a4paper, 12pt]{article}
\usepackage{amsfonts}
\usepackage{fancyhdr}
\usepackage{amssymb}
 \usepackage{mathrsfs,amsfonts,amsmath}
 \usepackage{color}

\makeatletter

\@addtoreset{equation}{section}
\makeatother
\newtheorem{Theorem}{Theorem}[section]
\newtheorem{Remark}[Theorem]{Remark}
\newtheorem{Lemma}[Theorem]{Lemma}
\newtheorem{Corollary}[Theorem]{Corollary}

\begin{document}
\title{\textbf{Exponential stability of the exact solutions and $\theta$-EM approximations
to neutral SDDEs with Markov switching}}
\author{ Guangqiang Lan$^a$\footnote{Email: langq@mail.buct.edu.cn. Supported by
National Natural Science Foundation of China (NSFC11026142) and
Beijing Higher Education Young Elite Teacher Project (BJYC34).}\quad
and\quad Chenggui Yuan$^b$
\\ \small $^a$School of Science, Beijing University of Chemical Technology, Beijing 100029, China
\\ \small $^b$Department of Mathematics, Swansea University, Singleton Park, Swansea, SA2 8PP, UK}

\date{}

\maketitle

\begin{abstract}
Exponential stability of the exact solutions as well as $\theta$-EM
($\frac{1}{2}<\theta\le 1$) approximations to neutral stochastic
differential delay equations with Markov switching will be
investigated in this paper. Sufficient conditions are obtained to
ensure the $p$-th moment ($p\ge1$) and almost sure exponential
stability of the exact solutions as well as $\theta$-EM
approximations ($p=2$). An example will be presented to support our
conclusions.
\end{abstract}



\section{Introduction}

\noindent

Many dynamical systems depend not only on present states but also on
past states. In such cases, stochastic differential delay equations
provide an important tool for describing such systems. More
generally, neutral stochastic  differential delay equations could be
used to describe the cases in which the delay argument occurs in the
derivative of the state variable. On the other hand, many practical
systems may experience abrupt changes in their structure and
parameters caused by phenomena such as component failures or
repairs, changing subsystem interconnections, and abrupt
environmental disturbances. The hybrid systems driven by continuous
time Markov chains have been developed to cope with such situations.
Motivated by hybrid systems, Kolmanovskii et al. \cite{KKMM}
introduced the neutral stochastic differential delay equations with
Markov switching (NSDDEswMS for short).

Let $(\Omega,\mathscr{F},P)$ be a probability space endowed with a
complete filtration $(\mathscr{F}_t)_{t\geq 0}$. Let
$d,m\in\mathbb{N}$ be arbitrarily fixed. For a given $\tau>0,$ let
$C([-\tau,0];\mathbb{R}^d)$ be a family of continuous functions
$\varphi$ from $[-\tau,0]$ to $\mathbb{R}^d,$ equipped with the
supremum norm
$||\varphi||=\sup_{\theta\in[-\tau,0]}|\varphi(\theta)|.$ Denote by
$C_{\mathscr{F}_0}^b([-\tau,0];\mathbb{R}^d)$ the family of bounded,
$\mathscr{F}_0$ measurable, $C([-\tau,0];\mathbb{R}^d)$-valued
random variables. If $x(t)$ is a continuous $\mathbb{R}^d$-valued
stochastic process on $t\ge 0,$ let $x_t=\{x(t+s):-\tau\le s\le0\}$
for $t\ge0$ which is regarded as a
$C([-\tau,0];\mathbb{R}^d)$-valued process. Let $r(t),t\ge0$ be a
right continuous Markov chain adapted to $(\mathscr{F}_t)_{t\geq 0}$
on the probability space taking values in a finite state space
$S=\{1,2,\cdots,N\}$ with generator $\Gamma=(\gamma_{ij})_{N\times
N}$ given by
$$P(r(t+\delta)=j|r(t)=i)=\gamma_{ij}\delta+o(\delta), i\neq j,$$
where $\delta>0.$ Here $\gamma_{ij}\ge0, i\neq j$ and $\gamma_{ii}:=\sum_{j\neq i}\gamma_{ij}.$

We consider the following NSDDEswMS
\begin{equation}\label{sde}d[X(t)-D(X(t-\tau),r(t))]
=f(X(t),X(t-\tau),t,r(t))dt+g(X(t),X(t-\tau),t,r(t))dB_t, \end{equation}
\noindent where the initial
$$x_0=\xi=\{\xi(\theta),\theta\in[-\tau,0]\}\in
C_{\mathscr{F}_0}^b([-\tau,0];\mathbb{R}^d),\quad r(0)=i_0\in S,$$
$(B_t)_{t\geq0}$ is an $m$-dimensional standard $\mathscr{F}_t$-Brownian motion independent
of $r(t)$, $D:(x,i)\in\mathbb{R}^d\times S\mapsto D(x,i)\in\mathbb{R}^d$,
$f:(x,y,t)\in\mathbb{R}^d\times\mathbb{R}^d\times[0,\infty)\mapsto f(x,y,t)\in\mathbb{R}^d$
and $g:(x,y,t)\in\mathbb{R}^d\times\mathbb{R}^d\times[0,\infty)\mapsto g(x,y,t)\in\mathbb{R}^d
\otimes\mathbb{R}^m$ are both Borel measurable functions.

It is well known that the Markov chain $r(t)$ can be represented as
a stochastic integral with respect to a Possion random measure:
$$dr(t)=\int_\mathbb{R}\bar{h}(r(t-),l)\nu(dt,dl),\quad t\ge0$$
with initial value $r(0)=i_0\in S,$ where $\nu(dt,dl)$ is a Possion
random measure with intensity $dt\times m(dl)$ in which $m$ is the
Lebesgue measure on $\mathbb{R}$ while the explicit definition of
$\bar{h}$ can be found in \cite{BBG,GAM}.

Kolmanovskii et al. \cite{KKMM} studied the existence and uniqueness
as well as the asymptotic moment boundedness and moment stability of
the solutions of the NSDDEswMS (\ref{sde}). In Mao et al.
\cite{MSY}, authors discussed the almost sure asymptotic stability
of the exact solutions of the NSDDEswMS  (\ref{sde}).

To make sure that equation (\ref{sde}) has a unique solution, which
is denoted by $X(t)\in\mathbb{R}^d,t\ge -\tau,$ \textit{throughout
this paper}, we assume that the coefficients satisfy local Lipschitz
condition, that is, for each $h>0,$ there is $L_h>0$ such that
\begin{equation}\label{local}|f(x,y,t,i)-f(\bar{x},\bar{y},t,i)|
\vee|g(x,y,t,i)-g(\bar{x},\bar{y},t,i)|\le L_h(|x-\bar{x}|+|y-\bar{y}|)\end{equation}
for all $(t,i)\in\mathbb{R}_+\times S$ and
$|x|\vee|y|\vee|\bar{x}|\vee|\bar{y}|\le h.$

We also assume in this paper that for all
$(x,y,i)\in\mathbb{R}^d\times\mathbb{R}^d\times S,$ there exists
$0<\beta<1$ such that
\begin{equation}\label{linear}|D(x,i)-D(y,i)|\le\beta|x-y|\end{equation}
and for all $(t,i)\in\mathbb{R}_+\times S,$
$$D(0,i)=0,\quad f(0,0,t,i)=0,\quad g(0,0,t,i)=0,$$
which implies that $x\equiv0$ is the trivial solution of equation
(\ref{sde}).

To study the exponential stability of equation (\ref{sde}), we
introduce the following operator $L$. If $V\in
C^{2,1}(\mathbb{R}^d\times\mathbb{R}_+\times S;\mathbb{R}_+),$
define $LV$ from
$\mathbb{R}^d\times\mathbb{R}^d\times\mathbb{R}_+\times S$ to
$\mathbb{R}_+$ by
\begin{equation}\label{LV}\aligned
LV(x,y,t,i):&=V_t(x-D(y,i),t,i)+V_x(x-D(y,i),t,i)f(x,y,t,i)\\&
\quad+\frac{1}{2}\textrm{trace}[
g^T(x,y,t,i)V_{xx}(x-D(y,i),t,j)g(x,y,t,i)]\\&
\quad+\sum_{j=1}^N\gamma_{ij}V(x-D(y,i),t,j),\endaligned\end{equation}
where
$$V_t(x,t,i)=\frac{\partial}{\partial t}V(x,t,i),\quad V_x(x,t,i)=
(\frac{\partial}{\partial x_1}V(x,t,i),\cdots,\frac{\partial}{\partial x_d}V(x,t,i))$$
and
$$V_{xx}(x,t,i)=(\frac{\partial^2}{\partial x_l\partial x_j}V(x,t,i))_{d\times d}.$$

Notice that if $V(x,t,i)=|x|^2,$ then
$$LV(x,y,t,i)=2\langle x-D(y,i),f(x,y,t,i)\rangle+||g(x,y,t,i)||^2,$$
where $|\cdot|$ is the Euclid norm in $\mathbb{R}^d$, $\langle
\cdot,\cdot\rangle$ denotes the inner product in $\mathbb{R}^d$ and
$||\cdot||$ stands for the Hilbert-Schmidt norm of a matrix.

When there are no neutral term and no delay, in \cite{Mao}, Mao
proved that $p$-th moment exponential stability plus linear growth
condition on coefficients $f$ and $g$ implies almost sure
exponential stability of the true solution of the corresponding
SDEs. In \cite{KKMM} Theorem 5.3, authors presented a sufficient
condition for almost sure exponential stability where the linear
growth condition is also needed. For other stability results of the
exact solutions, one can see e.g. \cite{LC,Mao1} and reference there
in. In this paper, we will first study the sufficient condition for
the almost sure exponential stability under the assumption of $p$-th
moment exponential stability, which is different with that of
\cite{KKMM}, where the linear condition for $g$ is needed (see the
following Remark \ref{rem}).

Since in general both the explicit solutions and the probability
distributions of the solutions are not known, discrete approximate
solutions become necessarily. A particularly important issue is to
determine the conditions under which the exact solution and
approximate solution share same stability properties. There are
plenty of literatures dealing with different types of numerical
approximations for SDEs, for example, \cite{BB,HMY,HMS,Mao1,WMS}
consider the moment or almost sure exponential stability of EM or
backward EM approximations, in \cite{LFM} the authors studied mean
square polynomial stability of EM and backward EM approximations,
while \cite{CW,Higham2,ZW} deal with the moment or almost sure
exponential stability of more general $\theta$-EM scheme. Recently,
in \cite{Lan}, the authors presented sufficient conditions for mean
square and almost sure exponential and polynomial stability of the
$\theta$-EM scheme. However, as far as we know, there are few
results on the exponential stability of the numerical approximations
of NSDDEswMS. So the other aim of this paper is to present
sufficient conditions for the exponential stability of numerical
approximations of NSDDEswMS.

Choose a step size $\Delta>0$ such that $\Delta=\tau/m$ for some
integer $m$. Define the corresponding discrete $\theta$
Euler-Maruyama ($\theta$-EM) approximation (or the so called
stochastic theta method) of the NSDDEswMS (\ref{sde}) as
\begin{equation}\label{SEM0} \aligned X_{j}=\xi(j\Delta), j=-m,-m+1,\cdots,0,\endaligned\end{equation}
while for $ k\ge0,$
\begin{equation}\label{SEM} \aligned X_{k+1}-D(X_{k+1-m},r((k+1)\Delta))&=X_k-D(X_{k-m},r(k\Delta))
+[(1-\theta)f(X_k,X_{k-m},k\Delta,r(k\Delta))\\&\quad+\theta f(X_{k+1},X_{k+1-m},(k+1)\Delta,r((k+1)\Delta))]\Delta\\&\quad
+g(X_k,X_{k-m},k\Delta,r(k\Delta))\Delta B_k,\endaligned\end{equation}

\noindent where $\theta\in [0,1]$ is a fixed parameter, $\Delta B_k:=B((k+1)\Delta t)-B(k\Delta t)$
is the increment of Brownian motion. Note that $\theta$-EM includes the classical EM method ($\theta=0$),
the backward EM method ($\theta=1$) and the so-called trapezoidal method ($\theta=\frac{1}{2}$).

Since the scheme (\ref{SEM}) is semi-implicit when $\theta>0$, to
make sure that the approximation scheme is well defined, a natural
assumption is that $f$ satisfies one-sided Lipschitz condition:
There exists $L>0$ such that for any
$x_1,x_2,y,\in\mathbb{R}^d,t\ge0,i\in S,$
\begin{equation}\label{Lip}
\langle x_1-x_2, f(x_1,y,t,i)-f(x_2,y,t,i)\rangle\le L|x_1-x_2|^2.
\end{equation}

Under this condition, if $L\theta\Delta<1,$ then the $\theta$-EM
scheme is well defined (see e.g. \cite{Lan,LFM,MS,WMS}).

We will investigate the sufficient conditions under which the
$\theta$-EM approximation is exponentially stable.

The rest of the paper is organized as follows. In Section 2, two
lemmas will be introduced to prove the following stability results.
In Section 3, sufficient conditions will be given to guarantee that
both $p$-th moment and almost sure exponential stability hold. In
Section 4, sufficient conditions will be presented to guarantee the
almost sure stability under the assumption that the exact solution
is exponentially stable in the $p$-th moment sense. Fianally, in
Section 5, we investigate the exponential stability of the
$\theta$-EM scheme. We show that when $\frac{1}{2}<\theta\le 1$,
both the mean square and almost sure exponential stability of the
$\theta$-EM scheme hold under given conditions. An example will be
presented to illustrate our theory.

\section{Preliminary}

Before we state our main results, let us present two useful lemmas
which will be used in the following Sections.

\begin{Lemma}\label{l1}
Let $p\ge1$. If (\ref{linear}) holds, then
$$|x-D(y,i)|^p\le (1+\beta)^{p-1}(|x|^p+\beta|y|^p),\quad\forall
(x,y,i)\in\mathbb{R}^d\times\mathbb{R}^d\times S.$$\end{Lemma}

For the proof, please see \cite{Mao}.

The second lemma is an elementary inequality which will be used in
the following Sections. One can find this inequality in \cite{MI}.

\begin{Lemma}\label{l2}
$$(a+b)^p\le
(1+c)^{p-1}(a^p+c^{1-p}b^p)$$\end{Lemma} holds for all $
a,b>0,p\ge1,c>0.$

\section{Sufficient conditions for exponential stability}
In this Section, we will study the $p$-th moment exponential
stability of $X(t)$. To do this, let us consider the following
conditions:

There exist functions $V\in
C^{2,1}(\mathbb{R}^d\times\mathbb{R}_+\times S;\mathbb{R}_+),$ $U\in
C(\mathbb{R}^d\times[-\tau,\infty);\mathbb{R}_+),$
$\gamma:\mathbb{R}_+\mapsto\mathbb{R}_+$ and $\psi\in
L^1(\mathbb{R}_+,\mathbb{R}_+)$ and positive constants
$c_1,\alpha_1,\alpha_2,\lambda$ such that $\int_0^\infty e^{\lambda
t}\gamma(t)dt<\infty$,
\begin{equation}\label{c1'}
c_1|x|^p\le V(x,t,i)
\end{equation}
and
\begin{equation}\label{c2'}
LV(x,y,t,i)\le\gamma(t)+(\psi(t)-\lambda)
V(x-D(y,i),t,i)-\alpha_1U(x,t)+\alpha_2U(y,t-\tau)
\end{equation}
holds for $\alpha_1\ge\alpha_2e^{\lambda\tau}, \beta
e^{\lambda\tau}<1.$

In general, conditions (\ref{local}) and (\ref{linear}) will only
guarantee a unique maximal local solution to SDE (\ref{sde}) for any
given initial value $\xi$ and $i_0.$ However, the additional
conditions (\ref{c1'}) and (\ref{c2'}) will guarantee that this
local solution is in fact a unique global solution. That is, we have
the following
\begin{Theorem}\label{exist}
Assume that conditions (\ref{local}), (\ref{linear}), (\ref{c1'})
and (\ref{c2'}) hold. Then there exist a unique global solution of
SDE (\ref{sde}).
\end{Theorem}
The proof is similar to that of Theorem A.1 in \cite{MSY}, so we
omit it here.

Now let us consider the exponential stability of equation
(\ref{sde}). We have the following
\begin{Theorem}\label{exact}
Assume that $p\ge1$, and conditions (\ref{local}), (\ref{linear}),
(\ref{c1'}) and (\ref{c2'}) hold. Then the exact solution of
equation (\ref{sde}) is $p$-th moment exponentially stable with
Lyapunov exponent no greater than $-\lambda$.
\end{Theorem}

\textbf{Proof} By the generalized It\^o's formula (see Skorohod
\cite{S} or Mao \cite{MY}), we have
$$\aligned d(e^{\lambda t}V(X(t)-D(X(t-\tau),r(t)),t,r(t)))&
=e^{\lambda t}[\lambda
V(X(t)-D(X(t-\tau),r(t)),t,r(t))\\&\quad+LV(X(t),X(t-\tau),t,r(t))]dt+dM_t.\endaligned$$
where
$$\aligned M_t&:=\int_0^te^{\lambda s}[V_x(X(s)-D(X(s-\tau),r(s)),s,r(s))g(X(s),X(s-\tau),s,r(s))dB_s\\&
\quad+\int_0^t\int_\mathbb{R}(V(X(s)-D(X(s-\tau),r(s)),s,r(s)+\bar{h}(r(s),l))
\\&\qquad-V(X(s)-D(X(s-\tau),r(s)),t,r(s)))\mu(ds,dl)]\endaligned$$
is a local martingale.

Then \begin{equation}\label{bds}\aligned &\quad e^{\lambda
t}V(X(t)-D(X(t-\tau),r(t)),t,r(t))\\& \le
V(X(0)-D(X(-\tau),i_0),0,i_0)+M_t\\&\quad+\int_0^te^{\lambda
s}[\gamma(s)+\psi(s)V(X(s)-D(X(s-\tau),r(s)),s,r(s))]ds\\& \le
V(X(0)-D(X(-\tau),i_0),0,i_0)+\int_0^te^{\lambda
s}\gamma(s)ds+M_t\\&\quad+\int_0^te^{\lambda
s}\psi(s)V(X(s)-D(X(s-\tau),r(s)),s,r(s))ds\\&\quad
-\alpha_1\int_0^te^{\lambda s}U(X(s),s)ds+\alpha_2\int_0^te^{\lambda
s}U(X(s-\tau),s-\tau)ds.\endaligned\end{equation}

Taking expectation on both sides (cutting by stopping time if
necessary). By conditions (\ref{c1'}) and (\ref{c2'}), we have
\begin{equation}\label{kongzhi}\aligned &\quad e^{\lambda
t}\mathbb{E}(V(X(t)-D(X(t-\tau),r(t)),t,r(t)))\\& \le
\mathbb{E}(V(X(0)-D(X(-\tau),i_0),0,i_0))+\int_0^te^{\lambda
s}\gamma(s)ds\\&\quad+\int_0^t\psi(s)e^{\lambda
s}\mathbb{E}(V(X(s)-D(X(s-\tau),r(s)),s,r(s)))ds\\&\quad
-\alpha_1\int_0^te^{\lambda
s}\mathbb{E}(U(X(s),s))ds+\alpha_2e^{\lambda\tau}\int_{-\tau}^0e^{\lambda
s}\mathbb{E}(U(X(s),s))ds\\&\quad+\alpha_2e^{\lambda\tau}\int_{0}^{t-\tau}e^{\lambda
s}\mathbb{E}(U(X(s),s))ds\\&\le C+\int_0^t\psi(s)e^{\lambda
s}\mathbb{E}(V(X(s)-D(X(s-\tau),r(s)),s,r(s)))ds,\endaligned\end{equation}
where
$$\aligned C&=\mathbb{E}(V(X(0)-D(X(-\tau),i_0),0,i_0))\\&\quad+\int_0^\infty e^{\lambda
s}\gamma(s)ds+\alpha_2e^{\lambda\tau}\int_{-\tau}^0e^{\lambda
s}\mathbb{E}(U(X(s),s))ds<\infty.\endaligned$$

We have used the fact that $\alpha_1\ge\alpha_2e^{\lambda\tau}$ in
the last inequality of (\ref{kongzhi}).

Now by Gronwall's lemma, we have
\begin{equation}\label{gronwall1}e^{\lambda
t}\mathbb{E}(V(X(t)-D(X(t-\tau),r(t)),t,r(t)))\le
Ce^{\int_0^t\psi(s)ds}\le
Ce^{\int_0^\infty\psi(s)ds}=:C'<\infty.\end{equation}

On the other hand, according to Lemma \ref{l2}, we have
$$\aligned e^{\lambda t}|X(t)|^p&\le e^{\lambda t}[(1-\beta)^{1-p}|X(t)-D(X(t-\tau),r(t))|^p
+\beta^{1-p}|D(X(t-\tau),r(t))|^p]\\&
\le\frac{(1-\beta)^{1-p}}{c_1}e^{\lambda
t}V(X(t)-D(X(t-\tau),r(t)),t,r(t))\\&\quad+\beta e^{\lambda
t}|X(t-\tau)|^p.\endaligned$$

Thus, for any $T>0,$ we have
$$\aligned \sup_{t\in[0,T]}e^{\lambda t}\mathbb{E}|X(t)|^p\le\frac{C'(1-\beta)^{1-p}}{c_1}
+\beta e^{\lambda \tau}[\sup_{t\in[-\tau,0]}e^{\lambda
t}\mathbb{E}|X(t)|^p+\sup_{t\in[0,T]}e^{\lambda
t}\mathbb{E}|X(t)|^p].\endaligned$$

Then
$$\aligned \sup_{t\in[0,T]}e^{\lambda t}\mathbb{E}|X(t)|^p\le(1-\beta e^{\lambda\tau})^{-1}[\frac{C'(1-\beta)^{1-p}}{c_1}
+\beta e^{\lambda \tau}\sup_{t\in[-\tau,0]}e^{\lambda
t}\mathbb{E}|X(t)|^p]<\infty,\endaligned$$ which implies
$$\limsup_{t\rightarrow\infty}\frac{\log\mathbb{E}|X(t)|^p}{t}\le-\lambda.\quad \square$$

\section{Almost sure exponential stability under the assumption of $p$-th moment exponential stability}

In this Section, we investigate the almost sure stability of exact
solution under the assumption that $p$-th moment exponential
stability holds for the exact solution.

Consider the following conditions:

There exist functions $V\in
C^{2,1}(\mathbb{R}^d\times\mathbb{R}_+\times S;\mathbb{R}_+)$,
$\varphi\in L^1(\mathbb{R}_+,\mathbb{R}_+)$ and
$h:\mathbb{R}_+\mapsto\mathbb{R}_+$ and positive constants
$c_1,c_2,\lambda,C_0,\alpha$ such that such that $\int_0^\infty
e^{\lambda t}h(t)dt<\infty$,
\begin{equation}\label{c1}
c_1|x|^p\le V(x,t,i)\le c_2|x|^p,
\end{equation}
\begin{equation}\label{c2}
LV(x,y,t,i)\le h(t)+(C_0+\varphi(t))V(x-D(y,i),t,i)
\end{equation}
and
\begin{equation}\label{c3}\aligned
&\quad|g^T(x,y,t,i)V_x(x-D(y,i),t,i)|^2\\&\quad
+\int_\mathbb{R}|V(x-D(y,i),t,i+\bar{h}(i,l))-V(x-D(y,i),t,i)|^2m(dl)\\&\le
\alpha V^2(x-D(y,i),t,i). \endaligned\end{equation}

We have the following
\begin{Theorem}\label{as}
Assume that the exact solution of equation (\ref{sde}) is $p$-th
moment exponentially stable ($p\ge1$) with Lyapunov exponent no
greater than $-\lambda$. If conditions (\ref{c1}), (\ref{c2}) and
(\ref{c3}) hold, and $\beta< e^{-\lambda\tau}$, then the exact
solution of equation (\ref{sde}) is almost surely exponentially
stable with Lyapunov exponent no greater than $-\frac{\lambda}{p}$.
\end{Theorem}

\textbf{Proof} By the generalized It\^o's formula, we have
$$\aligned &\quad dV(X(t)-D(X(t-\tau),r(t)),t,r(t))\\&=LV(X(t),X(t-\tau),t,r(t))dt\\&
\quad+V_x(X(t)-D(X(t-\tau),r(t)),t,r(t))g(X(t),X(t-\tau),t,r(t))dB_t\\&
\quad+\int_\mathbb{R}(V(X(t)-D(X(t-\tau),r(t)),t,r(t)+\bar{h}(r(t),l))\\&
\qquad-V(X(t)-D(X(t-\tau),r(t)),t,r(t)))\mu(dt,dl).\endaligned$$

Thus, for any $t\in[(k-1)\delta,k\delta),\ \delta>0$, we have
$$\aligned &\quad V(X(t)-D(X(t-\tau),r(t)),t,r(t))\\&
=V(X((k-1)\delta)-D(X((k-1)\delta-\tau),r((k-1)\delta)),(k-1)\delta,r((k-1)\delta))\\&
\quad+\int_{(k-1)\delta}^t LV(X(s),X(s-\tau),s,r(s))ds
+M_{1,t}+M_{2,t},\endaligned$$ where
$$M_{1,t}:=\int_{(k-1)\delta}^tV_x(X(s)-D(X(s-\tau),r(s)),s,r(s))g(X(s),X(s-\tau),s,r(s))dB_s$$
and
$$\aligned M_{2,t}&:=\int_{(k-1)\delta}^t\int_\mathbb{R}(V(X(s)-D(X(s-\tau),r(s)),s,r(s)+\bar{h}(r(s),l))\\&
\qquad-V(X(s)-D(X(s-\tau),r(s)),s,r(s)))\mu(ds,dl),\endaligned$$
$\mu(ds,dl)=\nu(ds,dl)-m(dl)ds$ is a martingale measure.

It is clear that both $M_{1,t}$ and $M_{2,t}$ are martingales on the
interval $[(k-1)\delta,k\delta).$ Denote
$M^*_{j,t}=\sup_{s\in[(k-1)\delta,t)} M_{j,s},\ j=1,2.$ Then by BDG
inequality, there exists $C$ (independent of $t, k$ and $\delta$)
such that
$$\aligned\mathbb{E}(M^*_{1,t})&\le C\mathbb{E}(\langle
M_{2,\cdot}\rangle_t^\frac{1}{2})\\&=
C\mathbb{E}[\int_{(k-1)\delta}^t|g^T(X(s),X(s-\tau),s,r(s))V_x(X(s)-D(X(s-\tau),r(s)),s,r(s))|^2ds]^\frac{1}{2}\\&\le
C\alpha\sqrt{\delta}\mathbb{E}\sup_{s\in[(k-1)\delta,t)}V(X(s)-D(X(s-\tau),r(s)),s,r(s)).\endaligned$$

Similarly, we have
$$\aligned\mathbb{E}(M^*_{2,t})&\le
C\alpha\sqrt{\delta}\mathbb{E}\sup_{s\in[(k-1)\delta,t)}V(X(s)-D(X(s-\tau),r(s)),s,r(s)).\endaligned$$

We have used the fact that
$$\aligned\langle
M_{2,\cdot}\rangle_t&=\int_{(k-1)\delta}^tds\int_\mathbb{R}|V(X(s)-D(X(s-\tau),r(s)),s,r(s)+\bar{h}(r(s),l))\\&
\qquad-V(X(s)-D(X(s-\tau),r(s)),s,r(s))|^2m(dl)\endaligned$$ in the
above inequality (see Situ \cite{situ} Lemma 62).

Then we have
$$\aligned&\quad\mathbb{E}\sup_{s\in[(k-1)\delta,t)}V(X(s)-D(X(s-\tau),r(s)),s,r(s))\\&\le
\mathbb{E}V(X((k-1)\delta)-D(X((k-1)\delta-\tau),r((k-1)\delta)),(k-1)\delta,r((k-1)\delta))\\&
\quad+\mathbb{E}\int_{(k-1)\delta}^t[h(s)+(C_0+\varphi(s))V(X(s)-D(X(s-\tau),r(s)),s,r(s))]ds\\&
\quad+2C\alpha\sqrt{\delta}\mathbb{E}\sup_{s\in[(k-1)\delta,t)}V(X(s)-D(X(s-\tau),r(s)),s,r(s))\\&
\le\mathbb{E}V(X((k-1)\delta)-D(X((k-1)\delta-\tau),r((k-1)\delta)),(k-1)\delta,r((k-1)\delta))\\&
\quad+
(C_0\delta+\int_{(k-1)\delta}^{k\delta}\varphi(s)ds+2C\alpha\sqrt{\delta})\mathbb{E}\sup_{s\in[(k-1)\delta,t)}
V(X(s)-D(X(s-\tau),r(s)),s,r(s))\\&\quad+\int_{(k-1)\delta}^th(s)ds.\endaligned$$

We can choose $\delta$ small enough such that
$\tilde{C}:=C_0\delta+\int_{(k-1)\delta}^{k\delta}\varphi(s)ds+2C\alpha\sqrt{\delta}<1.$

Thus,
$$\aligned&\quad\mathbb{E}\sup_{s\in[(k-1)\delta,t)}V(X(s)-D(X(s-\tau),r(s)),s,r(s))\\&\le
(1-\tilde{C})^{-1}[\mathbb{E}V(X((k-1)\delta)-D(X((k-1)\delta-\tau),r((k-1)\delta)),(k-1)\delta,r((k-1)\delta))\\&
\quad+\int_{(k-1)\delta}^th(s)ds]\\&\le
(1-\tilde{C})^{-1}[c_2\mathbb{E}|X((k-1)\delta)-D(X((k-1)\delta-\tau),r((k-1)\delta))|^p+\int_{(k-1)\delta}^th(s)ds].\endaligned$$

It follows that
$$\aligned&\quad\mathbb{E}\sup_{s\in[(k-1)\delta,t)}|X(s)-D(X(s-\tau),r(s))|^p\\&\le
c_1^{-1}(1-\tilde{C})^{-1}[c_2\mathbb{E}|X((k-1)\delta)
-D(X((k-1)\delta-\tau),r((k-1)\delta))|^p+\int_{(k-1)\delta}^th(s)ds].\endaligned$$

Since $X(t)$ is $p$-th moment exponentially stable, then by Lemma
\ref{l1}, for any $\varepsilon>0$, there exists $k_0$ large enough
such that for any $k\ge k_0$, we have
$$\aligned&\quad\mathbb{E}|X((k-1)\delta)-D(X((k-1)\delta-\tau),r((k-1)\delta))|^p\\&
\le
(1+\beta)^{p-1}(\mathbb{E}|X((k-1)\delta)|^p+\beta\mathbb{E}|X((k-1)\delta-\tau)|^p)\\&
\le(1+\beta)^{p-1}(e^{(-\lambda+\varepsilon)(k-1)\delta}+\beta
e^{(-\lambda+\varepsilon)[(k-1)\delta-\tau]})\\& \le
(1+\beta)^{p-1}(1+\beta
e^{\lambda\tau})e^{(-\lambda+\varepsilon)(k-1)\delta}.\endaligned$$

On the other hand, choosing $c=\frac{\beta}{1-\beta}$ in Lemma
\ref{l2}, we have
$$|X(t)|^p\le
(1-\beta)^{1-p}|X(t)-D(X(t-\tau),r(t))|^p+\beta^{1-p}|D(X(t-\tau),r(t))|^p$$
and $\int_{(k-1)\delta}^th(s)ds\le
e^{-\lambda(k-1)\delta}\int_0^\infty e^{\lambda s}h(s)ds,$ then
$$\aligned\mathbb{E}\sup_{s\in[(k-1)\delta,t)}|X(s)|^p&\le (1-\beta)^{1-p}c_1^{-1}(1-\tilde{C})^{-1}
[c_2(1+\beta)^{p-1}(1+\beta
e^{\lambda\tau})e^{(-\lambda+\varepsilon)(k-1)\delta}\\&
\quad+\int_{(k-1)\delta}^th(s)ds]+\beta\mathbb{E}\sup_{s\in[(k-1)\delta,t)}|X(s-\tau)|^p\\&
\le
C'e^{(-\lambda+\varepsilon)(k-1)\delta}+\beta\mathbb{E}\sup_{s\in[(k-1)\delta-\tau,t-\tau)}|X(s)|^p,\endaligned$$
where $C'=(1-\beta)^{1-p}c_1^{-1}(1-\tilde{C})^{-1}
[c_2(1+\beta)^{p-1}(1+\beta e^{\lambda\tau})+\int_0^\infty
e^{\lambda s}h(s)ds]<\infty$ is a constant independent of $k$ and
$t.$

Let $\tau/\delta=l$ be an integer. Then we have
$$\aligned a_{k,t}&\le b_k+\beta a_{k-l,t},\endaligned$$
where $a_{k,t}:=\mathbb{E}\sup_{s\in[(k-1)\delta,t)}|X(s)|^p,
b_k=C'e^{(-\lambda+\varepsilon)(k-1)\delta}.$

For any $t$, we can choose integers $n$ and $k$ such that
$[(k-1-nl)\delta,t-nl\delta]\subset [-\tau,0],$ which implies that
$n\ge\frac{(k-1)\delta}{\tau}.$ Notice also that $\beta
e^{(\lambda-\varepsilon)\tau}<1$ for any $\varepsilon>0.$ So we have
$$a_{k,t}\le\sum_{i=0}^{n-1}\beta^ib_{k-il}+\beta^{n}a_{k-nl,t}
\le\frac{C'e^{(-\lambda+\varepsilon)(k-1)\delta}}{1-\beta
e^{(\lambda-\varepsilon)\tau}}
+\beta^{\frac{(k-1)\delta}{\tau}}\sup|x_0|^p.$$

Therefore, for $k$ large enough,
$$\mathbb{E}\sup_{s\in[(k-1)\delta,t)}|X(s)|^p\le 2[\frac{C'e^{(-\lambda+\varepsilon)(k-1)\delta}}{1-\beta
e^{(\lambda-\varepsilon)\tau}}]\vee
[\beta^{\frac{(k-1)\delta}{\tau}}\sup|x_0|^p]\le
C{''}e^{-\lambda'(k-1)\delta},$$ where $C{''}:=\frac{2C'}{1-\beta
e^{(\lambda-\varepsilon)\tau}}\vee\sup|x_0|^p,$
$\lambda'=\lambda-\varepsilon.$

Hence, Chebyshev inequality yields that
$$P(\sup_{s\in[(k-1)\delta,t)}|X(s)|^p\ge
e^{-(\lambda-2\varepsilon)(k-1)\delta})\le\frac{\mathbb{E}
\sup_{s\in[(k-1)\delta,t)}|X(s)|^p}{e^{-(\lambda-2\varepsilon)(k-1)\delta}}\le
C{''}e^{-\varepsilon(k-1)\delta}.$$

By Borel-Cantelli lemma, we have for almost all $\omega\in\Omega,$
\begin{equation}\label{almost}\sup_{s\in[(k-1)\delta,k\delta)}|X(s)|
\le e^{-\frac{(\lambda-2\varepsilon)(k-1)\delta}{p}}\end{equation}
holds for all but finitely many $k.$ Hence there exists a
$k_0(\omega)$, for all $\omega\in\Omega$ excluding a P-null set, for
which equation (\ref{almost}) holds whenever $k\ge k_0$.
Consequently, for almost all $\omega\in\Omega$,
$$\frac{1}{t}\log|X(t)|\le-\frac{(\lambda-2\varepsilon)(k-1)\delta}{pt}
\le-\frac{(\lambda-2\varepsilon)(k-1)}{pk}$$
if $(k-1)\delta\le t\le k\delta$ and $k\ge k_0.$

This yields immediately that
$$\limsup_{t\rightarrow\infty}\frac{1}{t}\log|X(t)|
\le\limsup_{k\rightarrow\infty}-\frac{(\lambda-2\varepsilon)(k-1)}{pk}
=-\frac{(\lambda-2\varepsilon)}{p},\quad
a.s.$$

Letting $\varepsilon\rightarrow0.$ We complete the proof. $\square$

\begin{Remark}\label{rem}
Notice that our conditions (\ref{c2}) and (\ref{c3}) contain some
cases where $f$ and $g$ do not satisfy the linear growth condition.
For example, suppose $d=2, m=1, r>0.$ Let $D(y,i)$ satisfies
condition (\ref{linear}), and
$$g(x,y,t,i)=|x-D(y,i)|^r(-x_2+D_2(y,i),x_1-D_1(y,i))^T,$$
$$ f(x,y,t,i)=-|x-D(y,i)|^{2r}(x-D(y,i))^T-(x-D(y,i))^T,$$
where $D_j$ is the $j$-th coordinate of $D.$ It is clear that the
local Lipschitz condition (\ref{local}) holds for $f$ and $g$.
Moreover, by letting $V(x,t,i)=|x|^2$, we have
$$\aligned LV(x,y,t,i)&=-|x-D(y,i)|^{2r+2}-2|x-D(y,i)|^2\le-2|x-D(y,i)|^2\\&
\le-\frac{1}{1+c_0}|x|^2+\frac{\beta^2}{c_0}|y|^2,\quad \forall
c_0>\frac{\beta^2}{1-\beta^2}.\endaligned$$ Thus, by \cite{MSY}
Theorem A.1, there exist a unique global solution for equation
(\ref{sde}). On the other hand, according to \cite{KKMM} Theorem
5.2, equation (\ref{sde}) is exponentially stable in mean square
sense. Since our conditions (\ref{c1}), (\ref{c2}) and (\ref{c3})
hold for $p=2, V(x,t,i)=|x|^2$, then by Theorem \ref{as}, equation
(\ref{sde}) is also almost surely exponentially stable. However, the
linear growth condition for $g$ dose not hold in our case.
\end{Remark}

\begin{Remark}
By taking $U\equiv0$, condition (\ref{c2'}) becomes to
\begin{equation}\label{c22}
LV(x,y,t,i)\le\gamma(t)+(\psi(t)-\lambda) V(x-D(y,i),t,i).
\end{equation}
Obviously, it is stronger than condition (\ref{c2}).
\end{Remark}

Then by Theorem \ref{as} and Theorem \ref{exact}, we have the
following
\begin{Corollary}\label{cc1}
Suppose conditions (\ref{c1}), (\ref{c3}) and (\ref{c22}) hold. Then
the exact solution of equation (\ref{sde}) is $p$-th moment and
almost surely exponentially stable.
\end{Corollary}

\begin{Remark}
Notice that the convergence theorem of nonnegative semi-martingales
method could not be used here to prove the almost sure stability in
our case since we can not prove directly that in (\ref{bds}) the
integral term
$$\int_0^te^{\lambda
s}\psi(s)V(X(s)-D(X(s-\tau),r(s)),s,r(s))ds$$ is finite almost
surely. However, by using Theorem \ref{as}, we can overcome the
difficulty in proving the almost sure exponential stability in
Corollary \ref{cc1}.
\end{Remark}

\section{Exponential stability of $\theta$-EM solution (\ref{SEM0}) and (\ref{SEM})}

Let us now consider the exponential stability of the corresponding
$\theta$-EM approximation of the exact solution (\ref{sde}).

Since one-sided Lipschitz condition (\ref{Lip}) holds,
$$F(x,y,t,i):=x-D(y,i)-\theta\Delta f(x,y,t,i)$$
is well defined for all $0\le\theta\le1$. Then (\ref{SEM}) becomes
to \begin{equation}\label{dingyi}\aligned
F(X_{k+1},X_{k+1-m},(k+1)\Delta,r((k+1)\Delta))&=F(X_{k},X_{k-m},k\Delta,r(k\Delta))
\\&\quad+f(X_{k},X_{k-m},k\Delta,r(k\Delta))\Delta\\&\quad+g(X_{k},X_{k-m},k\Delta,r(k\Delta))\Delta
B_k.\endaligned\end{equation}

Denote $F_k:=F(X_{k},X_{k-m},k\Delta,r(k\Delta))$ for simplicity,
and $f_k,g_k$ in the following are defined in the same way.

Consider the following condition:
\begin{equation}\label{dandiao}
2\langle
x-D(y,i),f(x,y,t,i)\rangle+|g(x,y,t,i)|^2\le-C_1|x|^2+C_2|y|^2,
C_1>C_2>0.
\end{equation}

It is clear that under the Lipschitz conditions (\ref{local}) and
(\ref{dandiao}), there exists a unique global solution of equation
(\ref{sde}).

For the stability of $\theta$-EM numerical approximation
(\ref{SEM0}) and (\ref{SEM}), we have the following

\begin{Theorem}\label{numerical}
Assume that $\beta^2<\frac{1}{3+2\sqrt{2}}$, and conditions
(\ref{Lip}) and (\ref{dandiao}) hold for
$C_1>\frac{3+2\sqrt{2}}{1-(3+2\sqrt{2})\beta^2}C_2$. If
$\frac{1}{2}<\theta\le1$, then the $\theta$-EM numerical
approximation (\ref{SEM0}) and (\ref{SEM}) is mean square and almost
surely exponentially stable.
\end{Theorem}

\textbf{Proof} By (\ref{dingyi}), we have
$$\aligned|F_{k+1}|^2&=|F_k|^2+[2\langle X_k-D(X_{k-m},r(k\Delta)),f_k\rangle
+|g_k|^2+(1-2\theta)|f_k|^2\Delta]\Delta+M_k,\endaligned$$
where
\begin{equation}\label{M}\aligned M_k&:=|g_k\Delta B_k|^2
-|g_k|^2\Delta+2\langle f_k\Delta,g_k\Delta B_k\rangle+2\langle
F_k,g_k\Delta B_k\rangle.\endaligned\end{equation}

By condition (\ref{dandiao}), for any $c_0>0,$ if
$0<C<\frac{C_1}{1+c_0},$ then we claim that
\begin{equation}\label{budengshi}\aligned &\quad 2\langle
X_k-D(X_{k-m},r(k\Delta)),f_k\rangle
+|g_k|^2+(1-2\theta)|f_k|^2\Delta\\&
\le-C|F_k|^2+(\frac{C_1}{c_0}\beta^2+C_2)|X_{k-m}|^2.\endaligned\end{equation}

Indeed,
$$\aligned (2\theta-1)|f_k|^2\Delta
-C|F_k|^2&=[(2\theta-1)\Delta-C\theta^2\Delta^2]|f_k|^2\\&\quad+2C\theta\Delta\langle
X_k-D(X_{k-m},r(k\Delta)),f_k\rangle\\&\quad-C|X_k-D(X_{k-m},r(k\Delta))|^2\\&
=a|f_k+b(X_k-D(X_{k-m},r(k\Delta))|^2\\&\quad-(ab^2+C)|X_k-D(X_{k-m},r(k\Delta))|^2,\endaligned$$
where
$$a:=(2\theta-1)\Delta-C\theta^2\Delta^2,
\quad b:=\frac{C\theta\Delta}{a}.$$

Since $\frac{1}{2}<\theta\le1$ and $0<C<\frac{C_1}{1+c_0},$ we can
choose $\Delta$ small enough such that $a\ge 0$ and $ab^2+C\le
\frac{C_1}{1+c_0}$, and therefore
$$\aligned(2\theta-1)|f_k|^2\Delta-C|F_k|^2&\ge -\frac{C_1}{1+c_0}|X_k-D(X_{k-m},r(k\Delta))|^2\\&
\ge-\frac{C_1}{1+c_0}((1+c_0)|X_k|^2+\frac{1+c_0}{c_0}\beta^2|X_{k-m}|^2)\\&
=-C_1|X_k|^2-\frac{C_1}{c_0}\beta^2|X_{k-m}|^2\\& \ge 2\langle
X_k-D(X_{k-m},r(k\Delta)),f_k\rangle
+|g_k|^2\\&\quad-(\frac{C_1}{c_0}\beta^2+C_2)|X_{k-m}|^2.\endaligned$$

The second inequality holds because of Lemma \ref{l2}.

Then (\ref{budengshi}) holds for all $\frac{1}{2}<\theta\le1.$

So
$$|F_{k+1}|^2\le|F_k|^2-C\Delta|F_k|^2+(\frac{C_1}{c_0}\beta^2+C_2)\Delta|X_{k-m}|^2+M_k.$$

Then for any $A>1,$
$$\aligned&\quad A^{(k+1)\Delta}|F_{k+1}|^2-A^{k\Delta}|F_{k}|^2\\&
\le
A^{(k+1)\Delta}[|F_{k}|^2(1-C\Delta)+(\frac{C_1}{c_0}\beta^2+C_2)|X_{k-m}|^2+M_k]-A^{k\Delta}|F_{k}|^2\\&
=A^{(k+1)\Delta}|F_{k}|^2(1-C\Delta-A^{-\Delta})+(\frac{C_1}{c_0}\beta^2+C_2)\Delta
A^{(k+1)\Delta}|X_{k-m}|^2+A^{(k+1)\Delta}M_k.\endaligned$$

For simplicity, denote $R_1=1-C\Delta-A^{-\Delta},
R_2=\frac{C_1}{c_0}\beta^2+C_2.$ Then
$$\aligned A^{k\Delta}|F_{k}|^2&\le|F_0|^2+R_1\sum_{i=0}^{k-1}A^{(i+1)\Delta}|F_{i}|^2\\&
\quad+R_2\Delta\sum_{i=0}^{k-1}
A^{(i+1)\Delta}|X_{i-m}|^2+\sum_{i=0}^{k-1}
A^{(i+1)\Delta}M_i.\endaligned$$

Note that
$$\aligned\sum_{i=0}^{k-1}A^{(i+1)\Delta}|X_{i-m}|^2&=A^\tau\sum_{i=-m}^{-1}A^{(i+1)\Delta}|X_{i}|^2
+A^{\tau}\sum_{i=0}^{k-1}A^{(i+1)\Delta}|X_{i}|^2-A^\tau\sum_{i=k-m}^{k-1}A^{(i+1)\Delta}|X_{i}|^2\\&
\le A^\tau\sum_{i=-m}^{-1}A^{(i+1)\Delta}|X_{i}|^2
+A^{\tau}\sum_{i=0}^{k-1}A^{(i+1)\Delta}|X_{i}|^2.\endaligned$$

On the other hand, by the definition of $F_k$, condition
(\ref{dandiao}) and Lemma \ref{l2}, we know that
\begin{equation}\label{F}\aligned
|F_k|^2&\ge|X_k-D(X_{k-m},r(k\Delta))|^2+C_1\theta\Delta|X_k|^2-C_2\theta\Delta|X_{k-m}|^2\\&
\ge\frac{1}{1+c_0}|X_k|^2-\frac{\beta^2}{c_0}|X_{k-m}|^2+C_1\theta\Delta|X_k|^2-C_2\theta\Delta|X_{k-m}|^2\\&
=(\frac{1}{1+c_0}+C_1\theta\Delta)|X_k|^2-(\frac{\beta^2}{c_0}+C_2\theta\Delta)|X_{k-m}|^2.\endaligned\end{equation}

Note that $R_1|_{\Delta=0}=0$ and
$$R'_1(\Delta)=-C+A^{-\Delta}\log A\le0$$ for
$1<A\le e^{C}.$

Thus, $R_1\le0$ for $A$ close to 1.

Therefore,
$$\aligned A^{k\Delta}|F_{k}|^2&\le|F_0|^2+R_1
\sum_{i=0}^{k-1}A^{(i+1)\Delta}\Big[(\frac{1}{1+c_0}+C_1\theta\Delta)|X_i|^2-(\frac{\beta^2}{c_0}+C_2\theta\Delta)|X_{i-m}|^2)\Big]\\&
\quad+R_2\Delta
\sum_{i=0}^{k-1}A^{(i+1)\Delta}|X_{i-m}|^2+\sum_{i=0}^{k-1}
A^{(i+1)\Delta}M_i\\&
=|F_0|^2+K_1\sum_{i=0}^{k-1}A^{(i+1)\Delta}|X_{i}|^2+K_2\sum_{i=0}^{k-1}A^{(i+1)\Delta}|X_{i-m}|^2+\sum_{i=0}^{k-1}
A^{(i+1)\Delta}M_i\\&\le
|F_0|^2+K_2A^\tau\sum_{i=-m}^{-1}A^{(i+1)\Delta}|X_{i}|^2+(K_1+K_2A^\tau)\sum_{i=0}^{k-1}A^{(i+1)\Delta}|X_{i}|^2+\sum_{i=0}^{k-1}
A^{(i+1)\Delta}M_i,\endaligned$$ where
$$K_1=R_1(\frac{1}{1+c_0}+C_1\theta\Delta),$$
$$K_2=R_2\Delta-R_1(\frac{\beta^2}{c_0}+C_2\theta\Delta).$$

Now
$$\aligned K_1+K_2A^\tau&=R_1(\frac{1}{1+c_0}-A^\tau\frac{\beta^2}{c_0})+(R_1C_1\theta+R_2-R_1C_2\theta)\Delta\\&
=\Delta
[\frac{R_1}{\Delta}(\frac{1}{1+c_0}-A^\tau\frac{\beta^2}{c_0})+(R_2+R_1(C_1-C_2)\theta)].\endaligned$$

It is obvious that
$$\lim_{\Delta\rightarrow0}\frac{R_1}{\Delta}=\log A-C<0.$$

So if
$$h:=(\log A-C)(\frac{1}{1+c_0}-A^\tau\frac{\beta^2}{c_0})+R_2\le0,$$
then we have $K_1+K_2A^\tau\le0.$

On the other hand, since $\beta^2<\frac{1}{3+2\sqrt{2}},$ then there
exists $c_0>0$ such that
$$\beta^2<\frac{c_0}{c_0^2+3c_0+2}\le\frac{1}{3+2\sqrt{2}}.$$

Choose $c_0=\sqrt{2},$ then
$$c_0-(1+c_0)^2\beta^2-(1+c_0)\beta^2=\sqrt{2}(1-(3+2\sqrt{2})\beta^2)>0.$$

Since
$$C_1>\frac{3+2\sqrt{2}}{1-(3+2\sqrt{2})\beta^2}C_2=\frac{c_0(1+c_0)^2}{c_0-(1+c_0)^2\beta^2-(1+c_0)\beta^2}C_2,$$
then there exists $\varepsilon$ small enough such that
$$\frac{c_0(1-\varepsilon)-(1+c_0)^2\beta^2-
(1-\varepsilon)(1+c_0)\beta^2}{c_0(1+c_0)^2}C_1>C_2.$$

By taking $C=\frac{(1-\varepsilon)C_1}{1+c_0}$, we have
$$\frac{c_0C-(1+c_0)\beta^2C_1-(1+c_0)\beta^2C}{c_0(1+c_0)}>C_2.$$

That is,
$$\frac{C_1\beta^2}{c_0}+C_2-C(\frac{1}{1+c_0}-\frac{\beta^2}{c_0})<0.$$

Thus, $h\le0$ holds for some $A\ (>1)$ close to 1.

Then
$$\aligned A^{k\Delta}|F_{k}|^2&\le
|F_0|^2+K_2A^\tau\sum_{i=-m}^{-1}A^{(i+1)\Delta}|X_{i}|^2+\sum_{i=0}^{k-1}
A^{(i+1)\Delta}M_i.\endaligned$$

Take expectation on both sides, we have
$$A^{k\Delta}\mathbb{E}|F_{k}|^2\le
\mathbb{E}(|F_0|^2+K_2A^\tau\sum_{i=-m}^{-1}A^{(i+1)\Delta}|X_{i}|^2)=:K<\infty.$$

By (\ref{F}), we have
$$\aligned A^{k\Delta}\mathbb{E}|X_{k}|^2&\le
(\frac{1}{1+c_0}+C_1\theta\Delta)^{-1}
\Big(A^{k\Delta}\mathbb{E}|F_{k}|^2+(\frac{\beta^2}{c_0}+C_2\theta\Delta)A^{k\Delta}\mathbb{E}|X_{k-m}|^2\Big)\\&
\le K(1+c_0)+(1+c_0)(\beta^2/c_0+C_2\theta\Delta)A^\tau
A^{(k-m)\Delta}\mathbb{E}|X_{k-m}|^2.\endaligned$$

Notice that since $\beta^2\frac{1+c_0}{c_0}<\frac{1}{(1+c_0)^2}<1,$
we can choose $\Delta$ small enough and $A$ close to 1 such that
$q:=(1+c_0)(\beta^2/c_0+C_2\theta\Delta)A^\tau<1.$

Denote $a_k:=A^{k\Delta}\mathbb{E}|X_{k}|^2, b:=K(1+c_0).$

Then $$\aligned a_k&\le b+qa_{k-m}\\&\le
b\sum_{i=0}^{[k/m]}q^i+q^{[k/m]+1}a_{k-([k/m]+1)m}\\&\le
\frac{b}{1-q}+|\xi((k-([k/m]+1)m)\Delta)|^2.\endaligned$$

Then there exists $K'$ (independent of $k$) such that
\begin{equation}\label{mean}\aligned
A^{k\Delta}\mathbb{E}|X_{k}|^2\le
K'<\infty.\endaligned\end{equation}

This immediately yields that
$$\limsup_{k\rightarrow\infty}\frac{\log\mathbb{E}|X_{k}|^2}{k\Delta}\le-\log A<0.$$

Let us now consider the almost sure stability. By Chebyshev
inequality, inequality (\ref{mean}) implies that
$$P(|X_k|^2>e^{-k\Delta (\log A-\varepsilon}))\le K'e^{-k\Delta\varepsilon}, \forall k\ge 1, \varepsilon>0.$$

Then by Borel-Cantelli lemma, we see that for almost all
$\omega\in\Omega$
\begin{equation}\label{bound}|X_k|^2\le e^{-k\Delta (\log A-\varepsilon)}\end{equation}
holds for all but finitely many $k$. Thus, there exists a
$k_0(\omega),$ for all $\omega\in\Omega$ excluding a $P$-null set,
for which (\ref{bound}) holds whenever $k\ge k_0$.

Therefore, for almost all $\omega\in\Omega$,
\begin{equation}\label{bound1}\frac{1}{k\Delta}\log|X_k|\le-\frac{\log A-\varepsilon}{2}\end{equation}
whenever $k\ge k_0$. Letting $k\rightarrow\infty$ and
$\varepsilon\downarrow0$. Then
\begin{equation}\label{bound2}\limsup_{k\rightarrow\infty}\frac{1}{k\Delta}\log|X_k|\le-\frac{\log A}{2}.\end{equation}

The proof is then complete. $\square$

\begin{Remark}
When there is no Markovian switching, in \cite{WC}, the authors
considered the mean square stability of the semi implicit Euler
method for NSDDEs. The numerical method they considered is a special
case of our $\theta$-EM scheme (\ref{SEM}) ($\theta=1$). The
sufficient conditions are the global Lipschitz continuity of both
$f$ and $g$ plus the monotonicity condition (3.4) on $f$. However,
we only need the local Lipschitz condition on the coefficients $f$
and $g$, the one sided Lipschitz condition (\ref{dandiao}) on $f$
and the monotonicity condition (\ref{dandiao}). Moreover, they only
got the mean square stability while we get the mean square
exponential stability of $\theta$-EM scheme. We also remark that in
\cite{WC} condition ($C1$) implies the linear growth condition
($C2$) since $f(t,0,0)=g(t,0,0)\equiv0.$ That is, condition ($C2$)
is not necessary.
\end{Remark}

Let us now give an example to illustrate the theory.

\textbf{Example} Let $r(t)$ be a Markov Chain on the state space
$S=\{1,2\}$ with generator $\Gamma$. Consider the following scalar
neutral stochastic differential delay equation with Markov switching
\begin{equation}
\aligned d[X(t)-\frac{1}{6r(t)}\sin
X(t-\tau)]&=(-6X(t)-X^5(t)-\frac{1}{2}\sin
X(t-\tau))r(t)dt\\&\quad+2\frac{X^3(t)X(t-\tau)}{(1+X^2(t-\tau))r(t)}dB_t\endaligned
\end{equation}
with the initial value $x_0(\theta)=\xi(\theta)\equiv1,\theta\in
[-1,0].$ It is clear that the coefficients
$f(x,y,i)=-(6x+x^5+\frac{1}{2}\sin y)i$ and
$g(x,y,i)=2\frac{x^3y}{(1+y^2)i}$ satisfy the local Lipschitz
condition, and $D(x,i)=\frac{1}{6i}\sin x$ satisfies condition
(\ref{linear}) for $\beta_0=\frac{1}{6}$
($<\frac{1}{\sqrt{3+2\sqrt{2}}}$). Moreover, $f$ satisfies condition
(\ref{Lip}).

On the other hand, we have $|g(x,y,i)|^2\le|x|^6,$ and
$$\aligned2(x-D(y,i))f(x,y,i)&=-12ix^2-2ix^6+(2-i)x\sin y+\frac{1}{3}x^5\sin y+\frac{1}{6}\sin^2y\\&
\le-12ix^2-2ix^6+(2-i)(2x^2+\frac{\sin^2
y}{8})\\&\quad+\frac{1}{3}(\frac{5}{6}x^6+\frac{1}{6}\sin^6y)+\frac{1}{6}\sin^2y\\&
\le-(8+2i)x^2-(2i-\frac{5}{18})x^6+(\frac{2-i}{8}+\frac{1}{18}+\frac{1}{6})y^2.\endaligned$$

Then
$$2(x-D(y,i))f(x,y,i)+|g(x,y,i)|^2\le-10x^2+\frac{25}{72}y^2.$$

Thus (\ref{dandiao}) holds for $C_1=10,$ $C_2=\frac{25}{72}.$ By
\cite{MSY} Theorem A.1 (Take $\gamma\equiv0,
U(x,i)=\frac{25}{72}x^2$ there), there exists a unique global
solution $X(t)$ on $t\ge-1.$ By taking $\beta_0=\frac{1}{6}$ in
(\ref{linear}) and
$V(x,t,i)=U(x,i)=|x|^2,\gamma=\psi=0,\lambda=1,\alpha_1=8,\alpha_2=\frac{5}{12},
$ we know that (\ref{c2'}) holds. Then by Theorem \ref{exact}, the
exact solution is mean square exponentially stable. On the other
hand, since
$\frac{C_1}{C_2}=28.8>\frac{3+2\sqrt{2}}{1-(3+2\sqrt{2})\beta_0^2},$
then by Theorem \ref{numerical}, the $\theta$-EM approximation is
also mean square and almost surely exponentially stable.

However, the coefficient of the diffusion term does not satisfy the
linear growth condition $|g(x,y,i)|\le k(|x|+|y|)$ in this case.
That is, we can not get the mean stability from Theorem 3.1 in
\cite{WC}.

\end{document}